\documentclass[final,1p,times]{elsarticle}

\usepackage{amssymb}
 \usepackage{amsthm}
\usepackage{amscd}
\usepackage{amsmath}
\usepackage{amsfonts}
\usepackage{amssymb}
\usepackage{graphicx}
\newtheorem{theorem}{Theorem}
\newtheorem{proposition}[theorem]{Proposition}
\newtheorem{lemma}[theorem]{Lemma}
\newtheorem{corollary}[theorem]{Corollary}

\usepackage{mathrsfs}
\usepackage{titletoc}

\newcommand{\la}{\Delta}

\newcommand{\ra}{\rightarrow}
\newcommand{\p}{\partial}
\newcommand{\f}{\frac}

\newcommand{\be}{\begin{equation}}
\renewcommand{\ra}{\rightarrow}
\newcommand{\ee}{\end{equation}}
\newcommand{\bea}{\begin{eqnarray}}
\newcommand{\eea}{\end{eqnarray}}
\newcommand{\bna}{\begin{eqnarray*}}
\newcommand{\ena}{\end{eqnarray*}}

\renewcommand{\le}{\left}
\newcommand{\ri}{\right}

\journal{***}

\begin{document}

\begin{frontmatter}
\title{A heat flow for the mean field equation on a finite graph}

\author{Yong Lin}
\ead{yonglin@mail.tsinghua.edu.cn}
\address{Yau Mathematical Sciences Center, Tsinghua University, Beijing 100084, China}

\author{Yunyan Yang\footnote{Corresponding author}}
 \ead{yunyanyang@ruc.edu.cn}

\address{Department of Mathematics,
Renmin University of China, Beijing 100872, P. R. China}

\begin{abstract}
Inspired by works of Cast\'eras (Pacific J. Math., 2015), Li-Zhu (Calc. Var., 2019) and Sun-Zhu
(Calc. Var., 2020), we propose a heat flow for the mean field equation on a connected finite graph
$G=(V,E)$. Namely
$$\le\{\begin{array}{lll}
\p_t\phi(u)=\Delta u-Q+\rho \f{e^u}{\int_Ve^ud\mu}\\[1.5ex]
u(\cdot,0)=u_0,
\end{array}\ri.$$
where $\Delta$ is the standard graph Laplacian, $\rho$ is a real number, $Q:V\ra\mathbb{R}$ is a function
satisfying $\int_VQd\mu=\rho$, and  $\phi:\mathbb{R}\ra\mathbb{R}$ is one of certain smooth functions
including $\phi(s)=e^s$. We prove that for any initial data $u_0$ and any $\rho\in\mathbb{R}$, there exists a unique solution
$u:V\times[0,+\infty)\ra\mathbb{R}$ of the above heat flow; moreover, $u(x,t)$ converges to some function
$u_\infty:V\ra\mathbb{R}$ uniformly in $x\in V$ as $t\ra+\infty$, and $u_\infty$ is a solution of the mean field equation
$$\Delta u_\infty-Q+\rho\f{e^{u_\infty}}{\int_Ve^{u_\infty}d\mu}=0.$$
Though $G$ is a finite graph, this result is still unexpected, even in the special case $Q\equiv 0$. Our approach reads as follows: the short time existence of the heat flow follows from the ODE theory;
various integral estimates give its long time existence; moreover we establish a
Lojasiewicz-Simon type inequality and use it to conclude the convergence of the heat flow.

\end{abstract}

\begin{keyword}
 Heat flow on graph\sep the Lojasiewicz-Simon inequality\sep
 mean field equation

\MSC[2010] 35R02\sep 34B45

\end{keyword}

\end{frontmatter}

\titlecontents{section}[0mm]
                       {\vspace{.2\baselineskip}}
                       {\thecontentslabel~\hspace{.5em}}
                        {}
                        {\dotfill\contentspage[{\makebox[0pt][r]{\thecontentspage}}]}
\titlecontents{subsection}[3mm]
                       {\vspace{.2\baselineskip}}
                       {\thecontentslabel~\hspace{.5em}}
                        {}
                       {\dotfill\contentspage[{\makebox[0pt][r]{\thecontentspage}}]}

\setcounter{tocdepth}{2}


\section{Introduction}
Let us start with the mean field equation on a closed Riemann surface $(\Sigma,g)$, which says
\be\label{Mean-0}-\Delta_g u+Q=\rho\f{e^u}{\int_\Sigma e^udv_g},\ee
where $\rho\in \mathbb{R}$ is a number, $Q:\Sigma\ra \mathbb{R}$ is a smooth function with $\int_\Sigma Qdv_g=\rho$,
and $\Delta_g$ is the Laplacian operator with respect to the metric $g$. This equation arises in various topics such as
conformal geometry \cite{KZ}, statistical mechanics \cite{Caglioti}, and the abelian Chen-Simons-Higgs model \cite{Tarantello,Caffarelli,
Ding-2}.  The existence of solutions of (\ref{Mean-0}) has been extensively investigated for several decades.
Landmark achievements have been obtained for the case $\rho\not=8k\pi$, $k\in\mathbb{N}$, \cite{Brezis-Merle,Chen-Lin,Ding-3,Li1999,
Li-Shafrir,Malchiodi2008,Struwe-Tarantello,Djadli2008}, and for the case $\rho=8\pi$ \cite{Ding-1}.

In 2015, Cast\'eras \cite{Casteras1,Casteras2} proposed and studied the following parabolic equation
\be\label{flow-0}\le\{\begin{array}{lll}
\f{\p}{\p t}e^u=\Delta_g u-Q+\rho\f{e^u}{\int_\Sigma e^udv_g}\\[1.5ex]
u(x,0)=u_0(x),
\end{array}\ri.\ee
where $u_0\in C^{2,\alpha}(\Sigma)$, $0<\alpha<1$, is the initial data, $\Delta_g$, $Q$ and $\rho$ are described as in (\ref{Mean-0}).
It is a gradient flow for the energy functional $J_\rho: W^{1,2}(\Sigma,g)\ra\mathbb{R}$ defined by
$$J_\rho(u)=\f{1}{2}\int_\Sigma |\nabla_g u|^2dv_g+\int_\Sigma Q udv_g-\rho\log\int_\Sigma e^udv_g,$$
where $\nabla_g$ is the gradient operator with respect to the metric $g$. It was proved by Cast\'eras that for any
$\rho\not=8k\pi$, $k=1,2,\cdots$, there exists some initial data $u_0$ such that $u(\cdot,t)$ converges to a function $u_\infty$ in
$W^{2,2}(\Sigma)$, where $u_\infty$ is a solution of the mean field equation (\ref{Mean-0});
For $\rho=8\pi$, a sufficient condition for convergence of the flow (\ref{flow-0}) was given by Li-Zhu \cite{Li-Zhu}. This gives a new proof of the result of
Ding-Jost-Li-Wang \cite{Ding-1}, which was extended by Chen-Lin \cite{Chen-Lin2}
 to a general critical case,  and generalized by Yang-Zhu \cite{Yang-Zhu} to a non-negative prescribed function
case. Recently, using a more refined analysis, Sun-Zhu \cite{Sun-Zhu} studied a modified version of (\ref{flow-0}), i.e. the parabolic equation
\be\label{flow-1}\le\{\begin{array}{lll}
\f{\p}{\p t}e^u=\Delta_g u-\f{8\pi}{{\rm Area}(\Sigma)}+8\pi\f{he^u}{\int_\Sigma he^udv_g}\\[1.5ex]
u(x,0)=u_0(x),
\end{array}\ri.\ee
where $h(x)\geq 0$, $h\not\equiv 0$ on $\Sigma$, and ${\rm Area}(\Sigma)=\int_\Sigma dv_g$ denotes the area of $\Sigma$. Clearly this is another
method of proving the result in \cite{Yang-Zhu}.

In this paper, we are concerned with the mean field equation on a finite graph. Let us fix some notations. Assume $G=(V,E)$ is a finite graph, where $V$ denotes the vertex set and $E$ denotes the edge set. 
For any edge $xy\in E$, we assume that its weight $w_{xy}>0$ and that $w_{xy}=w_{yx}$.
Let $\mu:V\ra \mathbb{R}^+$ be a finite measure. For any function $u:V\ra \mathbb{R}$, the Laplacian
of $u$ is defined as
\be\label{lap}\Delta u(x)=\f{1}{\mu(x)}\sum_{y\sim x}w_{xy}(u(y)-u(x)),\ee
where $y\sim x$ means $xy\in E$, or $y$ and $x$ are adjacent.
The associated gradient form reads
$$\Gamma(u,v)(x)=\f{1}{2\mu(x)}\sum_{y\sim x}w_{xy}(u(y)-u(x))(v(y)-v(x)).$$
Write $\Gamma(u)=\Gamma(u,u)$. We denote the length of its gradient by
\be\label{gr}|\nabla u|(x)=\sqrt{\Gamma(u)(x)}=\le(\f{1}{2\mu(x)}\sum_{y\sim x}w_{xy}(u(y)-u(x))^2\ri)^{1/2}.\ee
For any function $g:V\ra\mathbb{R}$,  an integral of $g$ over $V$ is defined by
\be\label{int}\int_V gd\mu=\sum_{x\in V}\mu(x)g(x).\ee
Let $W^{1,2}(V)$ be a Sobolev space including all real functions $u$ with the norm
$$\|u\|_{W^{1,2}(V)}=\le(\int_V(|\nabla u|^2+u^2)d\mu\ri)^{1/2}.$$
As an analog of (\ref{Mean-0}), the mean field equation on the finite graph $G$ reads as
\be\label{Mean-1}-\Delta u+Q=\rho\f{e^u}{\int_V e^ud\mu},\ee
where $\rho$ is a real number, $Q:V\ra \mathbb{R}$ is a function with $\int_V Qd\mu=\rho$,
and $\Delta$ is the graph Laplacian with respect to the measure $\mu$ as in (\ref{lap}). The equation (\ref{Mean-1})
can be viewed as a discrete version of (\ref{Mean-0}). Let $\phi: \mathbb{R}\ra\mathbb{R}$ be a $C^1$ function.
We propose the following heat flow
\be\label{heat-flow}\le\{\begin{array}{lll}
\f{\p}{\p t}\phi(u)=\Delta u-Q+\rho\f{e^u}{\int_V e^ud\mu}\\[1.5ex]
u(x,0)=u_0(x),\,\,\, x\in V.
\end{array}\ri.
\ee
This is an analog of (\ref{flow-1}).
Obviously it is a gradient flow for the functional $J_\rho: W^{1,2}(V)\ra \mathbb{R}$, which is defined as
\be\label{functional}J_\rho(u)=\f{1}{2}\int_V|\nabla u|^2d\mu+\int_VQud\mu-\rho
\log\int_Ve^ud\mu,\ee
where the notations (\ref{gr}) and (\ref{int}) are used.
Our main result is stated as follows:

\begin{theorem}\label{thm1}
 Let $G=(V,E)$ be a connected finite graph. Suppose $\phi:\mathbb{R}\ra\mathbb{R}$ is a $C^1$ function  satisfying
 \be\label{phi}
 \lim_{s\ra-\infty}\phi(s)=0,\quad
 \phi^\prime(s)>0\,\,{\rm for\,\,all}\,\, s\in\mathbb{R},\quad
 \inf_{s\in[0,+\infty)}\phi^\prime(s)>0.
 \ee
 Let $\rho$ be
any real number, and $Q$ be any function with $\int_VQd\mu=\rho$.  Then for any initial function $u_0:V\ra\mathbb{R}$, we have the following
assertions:\\
$(i)$ there exists a unique solution
$u:V\times[0,\infty)\ra\mathbb{R}$ of the heat flow (\ref{heat-flow});\\
$(ii)$  there exists some function $u_\infty:V\ra\mathbb{R}$
such that $u(\cdot,t)$ converges to $u_\infty$ uniformly in $x\in V$ as $t\ra+\infty$;
moreover $u_\infty$ is a solution of the mean field equation (\ref{Mean-1}).
\end{theorem}

There are two interesting special cases of results in Theorem \ref{thm1} as follows:
\begin{corollary}\label{Cor}
Let $G=(V,E)$ and $\phi$ be as in Theorem \ref{thm1}.
If $\int_Vfd\mu=0$, then for any initial function $u_0:V\ra\mathbb{R}$, the heat flow
$$\le\{\begin{array}{lll}
\f{\p}{\p t}\phi(u)=\Delta u-f\\[1.5ex]
u(\cdot,0)=u_0
\end{array}\ri.$$
has a solution $u:V\times[0,+\infty)\ra\mathbb{R}$. Moreover, there exists some function $u^\ast$ such that
$u(\cdot,t)$ converges to $u^\ast$ as $t\ra+\infty$ uniformly in $x\in V$, and that $u^\ast$ satisfies
$$\le\{\begin{array}{lll}
\Delta u^\ast=f\\[1.5ex]
\int_V\phi{(u^\ast)}d\mu=\int_V\phi(u_0)d\mu.&&
\end{array}\ri.$$
\end{corollary}

\begin{corollary}\label{Cor-2}
Let $G=(V,E)$ and $\phi$ be as in Theorem \ref{thm1}.
Then for any initial function $u_0:V\ra\mathbb{R}$, the heat flow
 $$\le\{\begin{array}{lll}
\f{\p}{\p t}\phi(u)=\Delta u\\[1.5ex]
u(\cdot,0)=u_0
\end{array}\ri.$$
has a solution $u:V\times[0,+\infty)\ra\mathbb{R}$; moreover, as $t\ra+\infty$, $u(\cdot,t)$ converges to a constant $c$ uniformly
in $x\in V$, in particular
$$\phi(c)=\f{1}{|V|}\int_V\phi(u_0)d\mu.$$
\end{corollary}

Obviously there are infinitely many examples of $\phi$ in Theorem \ref{thm1}.  A typical example is
$$\phi(s)=\le\{\begin{array}{lll}
e^{\alpha s}+\beta s^p&{\rm when}& s>0\\[1.5ex]
e^{\alpha s}&{\rm when}& s\leq 0,
\end{array}\ri.$$
where $\alpha>0$, $\beta\geq 0$ and $p>1$ are constants. Another one says for any real number $a>1$,
$$\phi(s)=\le\{\begin{array}{lll}
s^2+(\log a)(s+\cos s-1)+1&{\rm when}& s>0\\[1.5ex]
a^{s}&{\rm when}& s\leq 0.
\end{array}\ri.$$

Though $G=(V,E)$ is a finite graph, the results in Theorem \ref{thm1} are quite unexpected, even in special cases
$\rho=0$ and $Q\equiv 0$ (Corollaries \ref{Cor} and \ref{Cor-2}).
As for its proof, we find a way of thinking from Simon \cite{Simon},  Jendoubi \cite{Jendoubi},
Cast\'eras \cite{Casteras1,Casteras2}, Li-Zhu \cite{Li-Zhu}, and Sun-Zhu \cite{Sun-Zhu}.
 Firstly, we use the ODE theory to conclude the short time existence of the heat flow (\ref{heat-flow}).
  Secondly, we obtain the global existence of the flow  through
estimating the uniform bound of $\|u(\cdot,t)\|_{W^{1,2}(V)}$  for all time $t$.
This allows us to select a sequence of times $(t_n)\ra+\infty$ such that $u(\cdot,t_n)$ converges to some function
$u_\infty$ uniformly in $V$, where $u_\infty$ is a solution of the mean field equation (\ref{Mean-1}).
Thirdly, we establish a Lojasiewicz-Simon type inequality along the heat flow by employing
an estimate due to Lojasiewicz (\cite{Lojasiewitz}, Theorem 4, page 88), namely

\begin{lemma}\label{Lojasie}{\rm (Lojasiewicz, 1963).}
 Let $\Gamma:\mathbb{R}^\ell\ra\mathbb{R}$ be an analytic function in a neighborhood of a point $\mathbf{a}\in\mathbb{R}^\ell$
 with $\nabla\Gamma(\mathbf{a})=\mathbf{0}\in\mathbb{R}^\ell$.
  Then there exist $\sigma>0$ and $0<\theta<1/2$ such that
  $$\|\nabla \Gamma(\mathbf{y})\|\geq |\Gamma(\mathbf{y})-\Gamma(a)|^{1-\theta},\quad
  \forall \mathbf{y}\in\mathbb{R}^\ell,\,\,\|\mathbf{y}-\mathbf{a}\|<\sigma,
  $$
  where $\nabla \Gamma(\mathbf{y})=({\p_{y_1}}\Gamma(\mathbf{y}),\cdots,\p_{y_\ell}\Gamma(\mathbf{y}))$, and $\|\cdot\|$ stands for
  the standard norm of $\mathbb{R}^\ell$.
\end{lemma}
\noindent Finally, we conclude the uniform convergence of $u(\cdot,t)$ to $u_\infty$ as $t\ra+\infty$ with the help of the above
Lojasiewicz-Simon type inequality. Since the graph $G$ is finite, this inequality seems much simpler than that of
\cite{Simon,Jendoubi,Casteras1,Casteras2,Li-Zhu,Sun-Zhu}. Moreover, in our case, all integral estimates
look very concise and very easy to understand.

Note that if the mean field equation (\ref{Mean-1}) has a solution, so does the Kazdan-Warner equation
\cite{GLY1}. For such kind of equations, see for examples \cite{Ge,Keller-Schwarz,Huang-Lin-Yau,Liu-Yang,Sun-Wang,Zhu}.
According to \cite{GLY1}, its solvability needs some assumptions. While Theorem \ref{thm1}
implies that for any real number $\rho$, (\ref{Mean-1}) is solvable. One may
ask whether or not these two conclusions are consistent.
Let us answer this question. Suppose that $u_\infty$ is given as in Theorem \ref{thm1}
and satisfies the mean field equation (\ref{Mean-1}).
 Clearly $v=u_\infty-\log\int_Ve^{u_\infty}d\mu$ is a solution of the Kazdan-Warner equation
 \be\label{K-W}\Delta v=Q-\rho e^v.\ee
 By the assumption in Theorem \ref{thm1}, $\int_VQd\mu=\rho$. If we assume $Q\equiv c$ for some constant $c$, then
 $\rho=c|V|$, where $|V|=\sum_{x\in V}\mu(x)$ denotes the volume of the graph.
 It follows from (\cite{GLY1}, Theorems 2-4) that if $c>0$ or $c<0$, then (\ref{K-W}) has a solution.
This implies that the results of Theorem \ref{thm1} and those of \cite{GLY1} are consistent and do not contradict each other.

For flow on infinite graph, partial existence results for the mean field equation are obtained by Ge-Jiang \cite{Ge-Jiang}.
Also there is a possibility of solving problems in \cite{GLY2,GLY3,Han-Shao-Zhao,Hou,Lin-Yang,Man,Zhang-Zhao} by a method of heat flow. Throughout this paper, we often denote various constants by the same $C$ from line to line,
even in the same line. The remaining part of this paper is arranged as follows: In Section 2, we prove the short time
existence of the heat flow; In Section 3, we show the heat flow exists for all time $t\in[0,+\infty)$; In Section 4, we
establish a Lojasiewicz-Simon type inequality and use it to prove the uniform convergence of the heat flow
as $t\ra+\infty$.
As a consequence, the proof of Theorem \ref{thm1} is finished.

 \section{Short time existence}\label{sec2}

 In this section, using the theory of ordinary differential equation, we shall prove that the solution of
 the heat flow (\ref{heat-flow}) exists on a short time interval. Also we shall give several properties of  the heat flow.

 Since $G=(V,E)$ is a finite graph, we assume with no loss of generality that
 $V=\{x_1,\cdots,x_\ell\}$ for some integer $\ell\geq 1$. Then any function
 $u:V\ra\mathbb{R}$ can be represented by
 $\mathbf{y}=(y_1,\cdots,y_\ell)\in\mathbb{R}^\ell$ with $y_j=u(x_j)$ for $1\leq j\leq \ell$;
 moreover, we denote
 \be\label{Mu}\mathcal{M}(u)=\Delta u-Q+\f{\rho e^u}{\int_Ve^{u}d\mu}\ee
 and a map $\mathcal{F}:\mathbb{R}^\ell\ra\mathbb{R}^\ell$ by
 $\mathcal{F}(\mathbf{y})=(f_1(\mathbf{y}),\cdots,f_\ell(\mathbf{y}))$, where
 $f_j(\mathbf{y})=\mathcal{M}(u)(x_j)$ for $1\leq j\leq\ell$.
 Then the equation  (\ref{heat-flow}) is equivalent to the ordinary differential system
 \be\label{ODEsystem}\le\{\begin{array}{lll}
 \f{d}{dt}\phi(y_1)&=&f_1(\mathbf{y})\\
 &\vdots&\\
 \f{d}{dt}\phi(y_\ell)&=&f_\ell(\mathbf{y})\\[1.2ex]
 \mathbf{y}(0)&=&\mathbf{y}_0,
 \end{array}\ri.\ee
 where $\mathbf{y}_0=(u_0(x_1),\cdots,u_0(x_\ell))$ is the initial data. For the map $\mathcal{F}$, we have the following:
 \begin{lemma}\label{analytic}
 The map $\mathcal{F}:\mathbb{R}^\ell\ra\mathbb{R}^\ell$ is analytic.
 \end{lemma}
 \proof At $\mathbf{y}=(u(x_1),\cdots,u(x_\ell))$, we write
 $$f_j(\mathbf{y})=\mathcal{M}(u)(x_j)=\f{1}{\mu(x_j)}\sum_{z\sim x_j}w_{zx_j}(u(z)-u(x_j))-
 Q(x_j)+\f{\rho e^{u(x_j)}}{\sum_{i=1}^\ell\mu(x_i)e^{u(x_i)}}.$$
 Replacing $(u(x_1),\cdots,u(x_\ell))$ by $\mathbf{y}$ on the righthand side of the above equality, we have that $f_j$ is analytic,
 $j=1,\cdots,\ell$. $\hfill\Box$\\

 On the short time existence of solutions of (\ref{heat-flow}), we obtain
 \begin{lemma}\label{short-time}
 There exists some constant $T^\ast>0$ such that (\ref{heat-flow}) has a solution $u:V\times [0,T^\ast]\ra\mathbb{R}$.
 \end{lemma}
 \proof By the short time existence theorem of the ordinary differential equation (\cite{ODE}, page 250),
 there exist some $T^\ast>0$ and a $C^1$ map $\mathbf{y}:[0,T^\ast]\ra \mathbb{R}^\ell$ such that
 $\mathbf{y}(t)$ satisfies (\ref{ODEsystem}). Define $u(x_j,t)=y_j(t)$ for $1\leq j\leq \ell$. Then $u:V\times [0,T^\ast]
 \ra\mathbb{R}$ is a solution of (\ref{heat-flow}). $\hfill\Box$\\

 For any $\rho\in\mathbb{R}$, let $J_\rho: W^{1,2}(V)\ra\mathbb{R}$ be a functional defined as in (\ref{functional}).
 One can easily see that (\ref{heat-flow}) is a negative gradient flow of $J_\rho$. In particular
 \be\label{dJ}\langle dJ_\rho(u(\cdot,t)),\phi\rangle=-\int_V\mathcal{M}(u(\cdot,t))\phi d\mu,\quad\forall\phi\in W^{1,2}(V).\ee
  Along the heat flow (\ref{heat-flow}), there are two important quantities: one is invariant, the other is monotone,
 namely

 \begin{lemma}\label{prop1}
 $(i)$ For all $t\in[0,T^\ast]$,  we have an invariant quantity
 $$\int_V\phi(u(\cdot,t))d\mu=\int_V\phi(u_0)d\mu.$$
 $(ii)$ $J_\rho(\cdot,t)$ is monotone with respect to $t$, in particular, if $0\leq t_1< t_2\leq T^\ast$, then
 $$J_\rho(u(\cdot,t_2))\leq J_\rho(u(\cdot,t_1)).$$
 \end{lemma}

 \proof
 Since $u(x,t)$ is a solution of (\ref{heat-flow}), we have by calculating
 \bna
 \f{d}{dt}\int_V\phi(u(\cdot,t))d\mu&=&\int_V\phi^\prime(u)u_td\mu\\
 &=&\int_V\le(\Delta u-Q+\f{\rho e^u}{\int_Ve^ud\mu}\ri)d\mu\\
 &=&0.
 \ena
 This immediately implies the assertion $(i)$.

 By the integration by parts,
 \bea\nonumber
 \f{d}{dt}J_\rho(u(\cdot,t))&=&\int_V\nabla u\nabla u_td\mu+\int_VQu_td\mu-
 \f{\rho }{\int_Ve^{u}d\mu}\int_Ve^{u}u_td\mu\\\nonumber
 &=&-\int_V\mathcal{M}(u)u_td\mu\\\label{deriv}
 &=&-\int_V\phi^\prime(u)u_t^2d\mu\leq 0,
 \eea
 since $\phi^\prime(s)>0$ for all $s\in\mathbb{R}$. Here we denote $u_t={\p}u/{\p t}$.
 This concludes the assertion $(ii)$.  $\hfill\Box$

 \section{Long time existence}
 In this section, we prove the long time existence of the heat flow (\ref{heat-flow}).
 By Lemma \ref{short-time}, there exists some $T^\ast>0$ such that (\ref{heat-flow}) has a solution $u:V\times [0,T^\ast]\ra
 \mathbb{R}$.
 Let \be\label{T}T=\sup\le\{T^\ast>0: u:V\times[0,T^\ast]\ra\mathbb{R}\,\, {\rm solves\,\,} (\ref{heat-flow})\ri\}.\ee
 Clearly, (\ref{heat-flow}) has a solution $u:V\times[0,T)\ra\mathbb{R}$.

 \begin{proposition}\label{prop2}
 Let $T$ be defined as in (\ref{T}). Then
 there exists some constant $C$ independent of $T$ such that
 $$\|u(\cdot,t)\|_{W^{1,2}(V)}\leq C,\quad\forall t\in[0,T).$$
 \end{proposition}

 \proof
   We divide the proof into several steps.

  {\bf Step 1}. {\it There exists some constant $C$ independent of $T$ such that for all $x\in V$ and $t\in[0,T)$,
  $$u(x,t)\leq C.$$}
  By (\ref{phi}), $\phi^\prime(s)>0$ for all $s\in\mathbb{R}$ and there exists some constant $a>0$ such that
  $\phi^\prime(s)\geq a>0$ for all $s\in[0,+\infty)$.
  There would hold $\phi(s)\geq as$ for all $s\in \mathbb{R}$. Indeed, the mean value theorem implies
  $$\phi(s)-\phi(0)=\phi^\prime(\xi)s,$$
  where $\xi$ lies between $0$ and $s$. Hence $\phi(s)\geq as$ for all $s\geq 0$ since $\phi(0)>0$. Obviously
  $\phi(s)\geq as$ for all $s< 0$ since $\phi(s)>0$ for all $s\in\mathbb{R}$.
  This together with $(i)$ of Lemma \ref{prop1} leads to
  \bna
  u(x,t)&\leq&\f{1}{a}\phi(u(x,t)) \\
  &\leq&\f{1}{a\min_{x\in V}\mu(x)}\int_V\phi(u(\cdot,t))d\mu\\
  &=&\f{1}{a\min_{x\in V}\mu(x)}\int_V\phi(u_0)d\mu\ena
  for all $x\in V$. This finishes the first step.\\

  {\bf Step 2}. {\it There exists a constant $C$ independent of $T$ such that for any $t\in [0,T)$, one finds a subset
  $A_t\subset V$ satisfying $\|u(\cdot,t)\|_{L^\infty(A_t)}\leq  C$ and $|A_t|\geq C^{-1}$.}\\

  For any $\epsilon>0$ and $t\in[0,T)$, we define a set
  $$V_{\epsilon,t}=\le\{x\in V: \phi(u(x,t))<\epsilon\ri\}.$$
  This together with $(i)$ of Lemma \ref{prop1} and Step 1 leads to
  \bea\nonumber
  \int_V\phi(u_0)d\mu&=&\int_V\phi(u(\cdot,t))d\mu\\\nonumber
  &=&\int_{V_{\epsilon,t}}\phi(u(\cdot,t))d\mu+\int_{V\setminus V_{\epsilon,t}}\phi(u(\cdot,t))d\mu\\
  &\leq&\epsilon|V|+\phi(C)|V\setminus V_{\epsilon,t}|.\label{11}
  \eea
  Taking $\epsilon=\epsilon_0=\f{1}{2|V|}\int_V\phi(u_0)d\mu$, we conclude from (\ref{11}) that
  \be\label{Vt}|V\setminus V_{\epsilon_0,t}|\geq \f{1}{2\phi(C)}\int_V\phi(u_0)d\mu.\ee
  Set $A_t=V\setminus V_{\epsilon_0,t}$. For any $x\in A_t$, there holds $\phi(u(x,t))\geq \epsilon_0$. Since $\phi(s)\ra 0$ as $s\ra-\infty$,
  we find some real number $b$ such that $\phi(b)=\epsilon_0$. It follows that $u(x,t)\geq b$ for all $x\in A_t$.
  This together with Step 1 leads to
  $$\|u(\cdot,t)\|_{L^\infty(A_t)}\leq C$$
  If $C$ is chosen larger but independent of $T$, then we have by (\ref{Vt}) that $|A_t|\geq C^{-1}$. \\

  {\bf Step 3}. {\it There exists a positive constant $C$ independent of $T$ such that for all $t\in[0,T)$,
  there holds
  $$\int_Vu^2(\cdot,t)d\mu\leq C\int_V|\nabla u(\cdot,t)|^2d\mu+C.$$}
  Recalling the definition of the first eigenvalue of the negative Laplacian, namely
  $$\lambda_1=\inf_{v\in W^{1,2}(V),\,\int_Vvd\mu=0,\,v\not\equiv 0}\f{\int_V|\nabla v|^2d\mu}{\int_Vv^2d\mu}>0,$$
  we obtain for any $v\in W^{1,2}(V)$,
  \bea\nonumber\int_Vv^2d\mu&=&\int_V(v-\overline{v})^2d\mu+\int_V\overline{v}^2d\mu\\&\leq& \f{1}{\lambda_1}\int_V|\nabla v|^2d\mu+\f{1}{|V|}\le(\int_Vvd\mu\ri)^2,\label{poincare}\eea
  where $\overline{v}=\f{1}{|V|}\int_Vvd\mu$.
  By Step 2, one calculates along the heat flow (\ref{heat-flow}),
  \bea\nonumber
  \f{1}{|V|}\le(\int_Vu(\cdot,t)d\mu\ri)^2&=&\f{1}{|V|}\le(\int_{A_t}u(\cdot,t)d\mu+\int_{V\setminus A_t}u(\cdot,t)d\mu\ri)^2\\
  \nonumber&=&\f{1}{|V|}\le(\int_{A_t}u(\cdot,t)d\mu\ri)^2+\f{1}{|V|}\le(\int_{V\setminus A_t}u(\cdot,t)d\mu\ri)^2\\\nonumber&&+\f{2}{|V|}\int_{A_t}u(\cdot,t)d\mu\int_{V\setminus A_t}u(\cdot,t)d\mu\\
  \nonumber&\leq&\f{C^2|A_t|^2}{|V|}+\f{1}{|V|}\le(\int_{V\setminus A_t}u(\cdot,t)d\mu\ri)^2\\\label{est-1}
  &&+\f{C^2|A_t|^2}{\epsilon |V|}+\f{\epsilon}{|V|}\le(\int_{V\setminus A_t}u(\cdot,t)d\mu\ri)^2,
  \eea
  where $\epsilon$ is a positive constant to be determined later. Using the H\"older inequality, one has
  $$\le(\int_{V\setminus A_t}u(\cdot,t)d\mu\ri)^2\leq |V\setminus A_t|\int_Vu^2(\cdot,t)d\mu.$$
  This together with (\ref{poincare}) and (\ref{est-1}) implies
  \be\int_Vu^2(\cdot,t)d\mu\leq \f{1}{\lambda_1}\int_V|\nabla u(\cdot,t)|^2d\mu+
  \f{(1+\epsilon)|V\setminus A_t|}{|V|}\int_Vu^2(\cdot,t)d\mu+{C^2|V|}\le(1+\f{1}{\epsilon}\ri).\label{est-2}\ee
  By Step 2, we have $|A_t|\geq C^{-1}$ with $0<C^{-1}<|V|$. Taking $\epsilon=(2C)^{-1}/(|V|-C^{-1})$
  in (\ref{est-2}), one
  gets $(1+\epsilon)|V\setminus A_t|/|V|\leq 1-(2C|V|)^{-1}$, and thus
  $$(2C|V|)^{-1}\int_Vu^2(\cdot,t)d\mu\leq \f{1}{\lambda_1}\int_V|\nabla u(\cdot,t)|^2d\mu+
  C.$$
  This completes the proof of Step 3.\\

  {\bf Step 4}. {\it There exists a constant $C$ independent of $T$ such that
  $\|u(\cdot,t)\|_{W^{1,2}(V)}\leq C$ for all $t\in[0,T)$.}\\

  It follows from the Poincar\'e inequality and the Young inequality that
  \be\label{ep-1}\int_V|u(\cdot,t)-\overline{u}(t)|d\mu\leq \epsilon\int_V|\nabla u(\cdot,t)|^2d\mu+C,\ee
  where $\epsilon>0$ is chosen later, $C$ is a constant depending on $\epsilon$, but independent of $T$.
  For any fixed $\rho\in\mathbb{R}$, in view of (\ref{functional}), we have by using (\ref{ep-1}),
  \bea\nonumber
  J_\rho(u(\cdot,t))&=&J_\rho(u(\cdot,t)-\overline{u}(t))\\\nonumber
  &=&\f{1}{2}\int_V|\nabla u(\cdot,t)|^2d\mu+\int_VQ(u(\cdot,t)-\overline{u}(t))d\mu\\\nonumber
  &&\quad-\rho\log\le(\int_Ve^{u(\cdot,t)-\overline{u}(t)}d\mu\ri)\\\label{tot}
  &\geq&\le(\f{1}{2}-\epsilon\ri)\int_V|\nabla u(\cdot,t)|^2d\mu-C-\rho\log \int_Ve^{u(\cdot,t)-\overline{u}(t)}d\mu.
  \eea
  As in Section \ref{sec2}, we write $V=\{x_1,\cdots,x_\ell\}$. Let $\theta_i=\mu(x_i)/|V|$ and $s_i=u(x_i,t)$, $1\leq i\leq \ell$.
  Obviously $0<\theta_i<1$ for any $i$ and $\sum_{i=1}^\ell\theta_i=1$. Since $e^s$ is convex in $s\in\mathbb{R}$, we have
  \bna
  \f{1}{|V|}\int_Ve^{u(\cdot,t)}d\mu
  &=&\sum_{i=1}^\ell\f{\mu(x_i)}{|V|}e^{u(x_i,t)}\\
  &=&\sum_{i=1}^\ell\theta_i e^{s_i}\\
  &\geq&e^{\sum_{i=1}^\ell\theta_is_i}\\
  &=&e^{\overline{u}(t)},
  \ena
  where $\overline{u}(t)=\f{1}{|V|}\int_Vu(\cdot,t)d\mu$. This immediately gives for $t\in[0,T)$,
  \be\label{et-01}\log\int_Ve^{u(\cdot,t)-\overline{u}(t)}d\mu\geq \log|V|.\ee
  According to the Trudinger-Moser embedding (\cite{GLY1}, Lemma 6), for any real number $\beta>0$,
  there exists some constant
  $C$ depending only on  $\beta$ and the Graph $G$ such that
  $$\int_Ve^{\beta\f{(u(\cdot,t)-\overline{u}(t))^2}
  {\|\nabla u(\cdot,t)\|_2^2}}d\mu\leq C.$$
  As a consequence,
  \bea\nonumber
  \log\int_Ve^{u(\cdot,t)-\overline{u}(t)}d\mu&\leq&\log\int_Ve^{\f{(u(\cdot,t)-\overline{u}(t))^2}
  {4\epsilon\|\nabla u(\cdot,t)\|_2^2}+\epsilon\|\nabla u(\cdot,t)\|_2^2}d\mu\\
  &\leq& \epsilon\int_V|\nabla u(\cdot,t)|^2d\mu+C\label{in-2}
  \eea
  for some constant $C$ depending on $\epsilon$ and the graph $G$. Combining (\ref{et-01}) and (\ref{in-2}),
  we have for any fixed real number $\rho$,
  \be\label{t0}\rho\log\int_Ve^{u(\cdot,t)-\overline{u}(t)}d\mu\leq |\rho|\epsilon\int_V|\nabla u(\cdot,t)|^2d\mu+C.\ee
  Inserting (\ref{t0}) into (\ref{tot}) and taking $\epsilon={1}/{(4+4|\rho|)}$, we conclude
  \be\label{J-lower}J_\rho(u(\cdot,t))\geq \f{1}{4}\int_V|\nabla u(\cdot,t)|^2d\mu-C.\ee
  This together with $(ii)$ of Lemma \ref{prop1} implies
  $$\int_V|\nabla u(\cdot,t)|^2d\mu\leq C, \quad\forall t\in[0,T).$$
  In view of Step 3, we complete the finial step and the proof of the lemma. $\hfill\Box$\\

  Now we are in a position to prove the long time existence of the heat flow (\ref{heat-flow}).

  \begin{proposition}\label{prop3}
  Let $T$ be given as in (\ref{T}). Then $T=+\infty$.
  \end{proposition}

  \proof Suppose $T<+\infty$. By Proposition \ref{prop2} and the short time existence theorem of the ordinary differential equation (\cite{ODE}, page 250), $u(\cdot,t)$ can be uniquely extended to a time interval $[0,T_2]$ for some $T_2>T$. This contradicts
  the definition of $T$. Therefore $T=+\infty$. $\hfill\Box$\\

  {\it Completion of the proof of $(i)$ of Theorem \ref{thm1}}. An immediate consequence of Proposition \ref{prop3}.
  $\hfill\Box$

  \section{Convergence of the heat flow}
  In this section, we shall prove $(ii)$ of Theorem \ref{thm1}.
  Since $V$ is finite, all norms of the function space $W^{1,2}(V)$ are equivalent. Then it follows from Proposition \ref{prop2}
  that there exists a constant $C$ such that for all $t\in [0,+\infty)$,
  \be\label{u-bd}\|u(\cdot,t)\|_{L^\infty(V)}\leq C.\ee
  In view of (\ref{deriv}) and (\ref{J-lower}), we have
  $$\int_0^{+\infty}\int_V\phi^\prime(u)u_t^2d\mu dt\leq J_\rho(u_0)+C.$$
  This together with the finiteness of $V$ and $\phi^\prime(s)>0$ for all $s\in\mathbb{R}$ implies that there exists an increasing sequence $t_n\ra+\infty$ such that for all $x\in V$,
  \be\label{ten-0}\le.\phi^\prime(u(x,t))u_t^2(x,t)\ri|_{t=t_n}\ra 0\quad{\rm as}\quad n\ra\infty.\ee
  By (\ref{u-bd}), since $\phi^\prime(s)>0$ for all $s\in\mathbb{R}$, we obtain
  \be\label{lower}0<\min_{s\in[-C,C]}\phi^\prime(s)\leq\phi^\prime(u(x,t_n))\leq
  \max_{s\in[-C,C]}\phi^\prime(s),\quad\forall n\geq 1.\ee
  Combining (\ref{ten-0}) and (\ref{lower}), we conclude that for all $x\in V$,
  \be\label{t-0}\le.\f{\p}{\p t}\phi(u(x,t))\ri|_{t=t_n}\ra 0\quad{\rm as}\quad n\ra\infty.\ee
  Moreover, up to a subsequence, we can find some function $u_\infty: V\ra\mathbb{R}$ such that
  $u(x,t_n)$ converges to $u_\infty(x)$ uniformly in $x\in V$ as $n\ra\infty$. This together with (\ref{heat-flow})
  and (\ref{t-0}) leads to
  \be\label{u-infty}\mathcal{M}(u_\infty)=\Delta u_\infty-Q+\f{\rho e^{u_\infty}}{\int_Ve^{u_\infty}d\mu}=0\quad{\rm on}
  \quad V.\ee
  In conclusion, we found an increasing sequence $(t_n)\ra+\infty$ such that $u(x,t_n)\ra u_\infty(x)$ uniformly in
  $x\in V$ as $n\ra\infty$, where $u_\infty$ is a solution of the mean field equation (\ref{u-infty}).
  Hereafter we further prove that along the heat flow, $u(\cdot,t)$ converges to $u_\infty$ as $t\ra+\infty$ uniformly
  on $V$. For this purpose, we need an estimate due to Lojasiewicz, namely
  Lemma \ref{Lojasie}. The power of Lemma \ref{Lojasie} is shown in the following finite dimensional Lojasiewicz-Simon inequality.

  \begin{proposition}\label{prop4} Let $\sigma>0$ and $0<\theta<1/2$ be given as in Lemma \ref{Lojasie}, $\mathcal{M}(u)$ be defined as in (\ref{Mu}),
  and $\ell$ be the number of points of $V$.
  Along the heat flow (\ref{heat-flow}), if $\|u(\cdot,t)-u_\infty\|_{L^\infty(V)}<\sigma/\sqrt{\ell}$ for some fixed $t$, then there exists some constant $C$
  independent of $t$ such that
  $$|J_\rho(u(\cdot,t))-J_\rho(u_\infty)|^{1-\theta}\leq C\|\mathcal{M}(u)(\cdot,t)\|_{L^2(V)}.$$
  \end{proposition}
  \proof Assume $\|u(\cdot,t)-u_\infty\|_{L^\infty(V)}<\sigma/\sqrt{\ell}$ for some fixed $t$. For the sake of clarity, we denote $\mathbf{y}=(y_1,\cdots,y_\ell)=(u(x_1,t),\cdots,u(x_\ell,t))$,
   $\mathbf{a}=(u_\infty(x_1),\cdots,u_\infty(x_\ell))$, $\Gamma(\mathbf{y})=J_\rho(u(\cdot,t))$ and $\Gamma(\mathbf{a})=
   J_\rho(u_\infty)$. Clearly the function $\Gamma:\mathbb{R}^\ell\ra\mathbb{R}$ is analytic due to Lemma \ref{analytic}, and  $$\|\mathbf{y}-\mathbf{a}\|=\sqrt{\sum_{i=1}^\ell(y_i-a_i)^2}\leq \sqrt{\ell}\max_{1\leq i\leq\ell}|y_i-a_i|<\sigma.$$
   For any $1\leq i\leq \ell$, we define a function $e_i:V\ra\mathbb{R}$ by
    $$e_i(x)=\le\{\begin{array}{lll}
    1,&{\rm if}& x=x_i\\[1.5ex]
    0,&{\rm if}& x\not=x_i.
    \end{array}\ri.$$
    Let $\mathbf{e}_i$ be a unit vector in $\mathbb{R}^\ell$, whose $i$-th component is $1$ and the rest are $0$.
   In view of (\ref{dJ}), one calculates the partial derivative of the analytic function $\Gamma(y)$ as follows.
   For any $1\leq i\leq\ell$,
   \bea\nonumber
   \p_{y^i}\Gamma(\mathbf{y})&=&\lim_{h\ra 0}\f{1}{h}\le(\Gamma(\mathbf{y}+h\mathbf{e}_i)-\Gamma(\mathbf{y})\ri)\\
   &=&\lim_{h\ra 0}\f{1}{h}\le(J_\rho(u(x,t)+he_i(x))-J_\rho(u(x,t))\ri)\nonumber\\
   &=&dJ_\rho(u(x,t))(e_i(x))\nonumber\\
   &=&\int_V\mathcal{M}(u)(x,t)e_i(x)d\mu.\label{d-gamma}
   \eea
   This together with the fact $\sum_{i=1}^\ell\int_Ve_i^2d\mu=\sum_{i=1}^\ell\mu(x_i)=|V|$ leads to
   \bea\|\nabla \Gamma(\mathbf{y})\|&=&\sqrt{\sum_{i=1}^\ell\le(\p_{y_i}\Gamma(y)\ri)^2}\nonumber\\
   &\leq& \sqrt{\le(\int_V\mathcal{M}(u)^2d\mu\ri)\sum_{i=1}^\ell\int_Ve_i^2d\mu}\nonumber\\
   &=& \sqrt{|V|}\,\|\mathcal{M}(u)(\cdot,t)\|_{L^2(V)}.\label{na-1}\eea
   Similar to (\ref{d-gamma}), we have by (\ref{u-infty}) that for all $1\leq i\leq\ell$,
   \be\label{nab-Ga}\p_{y_i}\Gamma(\mathbf{a})=\int_V\mathcal{M}(u_\infty(x))e_i(x)d\mu=0.\ee
   In view of the definition of $\Gamma$, (\ref{na-1}) and (\ref{nab-Ga}), we obtain by applying Lemma \ref{Lojasie} that
   \bna
   |J_\rho(u(\cdot,t))-J_\rho(u_\infty)|^{1-\theta}&=&|\Gamma(\mathbf{y})-\Gamma(\mathbf{a})|^{1-\theta}\\
   &\leq&\|\nabla \Gamma(y)\|\\
   &\leq&\sqrt{|V|}\,\|\mathcal{M}(u)(\cdot,t)\|_{L^2(V)}.
   \ena
   This ends the proof of the proposition. $\hfill\Box$\\

   Finally we prove the uniform convergence of the heat flow (\ref{heat-flow}), namely
   \begin{proposition}\label{uniform}
   Along the heat flow (\ref{heat-flow}), there holds
   \be\label{uni}\lim_{t\ra+\infty}\int_V|u(\cdot,t)-u_\infty|^2d\mu=0.\ee
   \end{proposition}
   \proof
   For the proof of this proposition,
   we modify an argument of Sun-Zhu (\cite{Sun-Zhu}, Section 5). Suppose that $(\ref{uni})$ does not hold. Then there exists some constant $\epsilon_0>0$ and
   an increasing sequence of numbers
   $(t_n^\ast)$ such that $t_n^\ast>t_n$ and
   \be\label{contr}\int_V|u(\cdot,t_n^\ast)-u_\infty|^2d\mu\geq 2\epsilon_0,\ee
   where $(t_n)$ is given by (\ref{ten-0}) and satisfies $u(\cdot,t_n)\ra u_\infty$ uniformly on $V$.
   Obviously
   $$\lim_{n\ra \infty}\int_V|u(\cdot,t_n)-u_\infty|^2d\mu=0.$$
   Thus there exists $n_1\in\mathbb{N}$ such that if $n\geq n_1$, then
   \be\label{lim}\int_V|u(\cdot,t_n)-u_\infty|^2d\mu<\epsilon_0.\ee
   We {\it claim} that $J_\rho(u(\cdot,t))> J_\rho(u_\infty)$ for all $t\in [0,+\infty)$. Indeed, we have by
   $(ii)$ of Lemma \ref{prop1} that $J_\rho(u(\cdot,t))$ is decreasing with respect to $t$, and in particular
   $J_\rho(u(\cdot,t))\geq J_\rho(u_\infty)$ for all $t\geq 0$. Suppose there exists some
   $\tilde{t}>0$ such that $J_\rho(u(\cdot,\tilde{t}))= J_\rho(u_\infty)$. Then $J_\rho(u(\cdot,{t}))\equiv J_\rho(u_\infty)$
   and thus $u_t\equiv 0$ on $V$ for all $t\in[\tilde{t},+\infty)$. Hence $u(x,t)\equiv u_\infty(x)$ for all $x\in V$ and all
   $t\in[\tilde{t},+\infty)$, which contradicts (\ref{contr}). This confirms our claim $J_\rho(u(\cdot,t))> J_\rho(u_\infty)$
   for all $t\geq 0$.

    For any $n\geq n_1$, we define
   $$s_n=\inf\le\{t>t_n: \|u(\cdot,t)-u_\infty\|_{L^2(V)}^2\geq 2\epsilon_0\ri\}.$$
   It follows from (\ref{contr}) that $s_n<+\infty$, and that for all $t\in[t_n,s_n)$,
   \be\label{equ-1}\int_V|u(\cdot,t)-u_\infty|^2d\mu<2\epsilon_0=\int_V|u(\cdot,s_n)-u_\infty|^2d\mu.\ee
   For $t\in[t_n,s_n)$, we calculate by (\ref{deriv}), Proposition \ref{prop4}, the fact $u_t=(\phi(u))^{-1}\mathcal{M}(u)$,
   and (\ref{u-bd}) that
   \bna
   -\f{d}{dt}(J_\rho(u(\cdot,t))-J_\rho(u_\infty))^{\theta}&=&-\theta
   (J_\rho(u(\cdot,t))-J_\rho(u_\infty))^{\theta-{1}}\f{d}{dt}J_\rho(u(\cdot,t))\\
   &=&\theta
   (J_\rho(u(\cdot,t))-J_\rho(u_\infty))^{\theta-{1}}\int_V\mathcal{M}(u)u_td\mu\\
   &\geq&C\f{\int_V(\phi^\prime(u))^{-1}\mathcal{M}^2(u)d\mu}{\|\mathcal{M}(u)\|_{L^2(V)}}\\
   &\geq& C\|u_t\|_{L^2(V)}.
   \ena
   Hence
   \be\label{est-4}\int_{t_n}^{s_n}\|u_t\|_{L^2(V)}dt\leq C(J_\rho(u(\cdot,t_n))-J_\rho(u_\infty))^{\theta}.\ee
   By the H\"older inequality,
   \be\label{est-5}\f{d}{dt}\le(\int_V|u(\cdot,t)-u_\infty|^2d\mu\ri)^{1/2}=\f{1}{\|u(\cdot,t)-u_\infty\|_{L^2(V)}}\int_V(u-u_\infty)u_td\mu\leq
   \le(\int_Vu_t^2d\mu\ri)^{1/2}.\ee
   Combining (\ref{est-4}) and (\ref{est-5}), we have
   $$\|u(\cdot,s_n)-u_\infty\|_{L^2(V)}-\|u(\cdot,t_n)-u_\infty\|_{L^2(V)}\leq C(J_\rho(u(\cdot,t_n))-J_\rho(u_\infty))^{\theta}.$$
   This together with (\ref{lim}) and (\ref{equ-1}) leads to
   $$\epsilon_0\leq C(J_\rho(u(\cdot,t_n))-J_\rho(u_\infty))^{\theta},$$
   which is impossible if $n$ is chosen sufficiently large,
   since $J_\rho(u(\cdot,t_n))\ra J_\rho(u_\infty)$ as $n\ra\infty$.
   This confirms (\ref{uni}). $\hfill\Box$\\

   {\it Completion of the proof of $(ii)$ of Theorem \ref{thm1}}. Recalling  $V=\{x_1,\cdots,x_\ell\}$, one concludes
   from Proposition \ref{uniform} that
      $$\lim_{t\ra+\infty}\sum_{i=1}^\ell \mu(x_i)|u(x_i,t)-u_\infty(x_i)|^2=0.$$
      Since for all $j\in\{1,\cdots,\ell\}$, there holds
      \bna
      |u(x_j,t)-u_\infty(x_j)|\leq \f{1}{\min_{x\in V}\mu(x)}\sum_{i=1}^\ell \mu(x_i)|u(x_i,t)-u_\infty(x_i)|^2,
      \ena
      one comes to a conclusion that $u(x,t)$ converges to $u_\infty(x)$ uniformly in $x\in V$ as $t\ra+\infty$.
   By (\ref{u-infty}), $u_\infty$ is a solution of (\ref{Mean-1}). Thus the proof of Theorem \ref{thm1} is completely
   finished. $\hfill\Box$\\

     {\bf Acknowledgements.} Yong Lin is partly supported by the National Science Foundation of China (Grant No. 12071245).
     Yunyan Yang is partly supported by the  National Science Foundation of China (Grant No. 11721101) and
     National Key Research and Development Project SQ2020YFA070080.
 Both of the two authors are supported by the National Science Foundation of China (Grant No. 11761131002).

\end{document}